\documentclass [10pt,a4paper,draft] {article}
\usepackage [applemac] {inputenc}
\usepackage[french] {babel}
\usepackage {amsfonts}
\usepackage {amsmath}
\usepackage {amssymb}

\begin{document}

\begin{center}
\begin{LARGE}
Note sur les invariants du groupe affine
\end{LARGE}
\bigskip

par : Mustapha RA\"IS  (Poitiers)
\end{center}

\bigskip
\vspace{10 mm}


\medskip

	Cette note m'a \'et\'e sugg\'er\'ee par :

\medskip

- Une proposition de E.M. Baruch (\cite{B}) : 

\medskip
	Soit $P$ le sous-groupe de $GL_n(K)\ (K = \Bbb{R}$ ou $\Bbb{C}$) constitu\'e par les matrices dont la
derni\`ere ligne est $(0,0,\ldots ,0,1)$. 

\bigskip
\noindent
\textbf{Key Proposition} (\cite{B}, page 209)~: Soit $T$ une distribution $P$-invariante sur l'ouvert des
\'el\'ements r\'eguliers de $GL_n(K)$. Alors $T$ est $GL_n(K)$-invariante.

\vskip 7mm
- Un lemme de T. Levasseur et J. T. Stafford (\cite{L-S}) : 

\medskip
	Soit $\mathfrak{g}$ une alg\`ebre de Lie (de dimension finie sur $K = \Bbb{R}$) qui admet une forme
bilin\'eaire $B$, invariante et non d\'eg\'en\'er\'ee. Lorsque $f : \mathfrak{g} \longrightarrow \Bbb{R}$ est
une fonction d\'erivable invariante (par le groupe adjoint de $\mathfrak{g}$), on sait que : $[\nabla
f(x),x]=0$ pour tout $x$ dans $\mathfrak{g}$, o\`u $\nabla f$ est le gradient de $f$, calcul\'e au moyen de
$B$ :
$$
	B(\nabla f(x),y) = \, <df(x),y>\, = \big( {d\over dt}\big)_0 \ f(x+ty)
$$
pour tout $y$ dans $\mathfrak{g}$.

\smallskip
	Soit $(e_i)_i$ une base de $\mathfrak{g}$, et soit $(e^*_i)_i$ la base duale au sens : $B(e_i,e^*_j) =
\delta_{ij}$. Alors : 
$$
	\nabla f(x) = \sum_i\ B(\nabla f(x), e^*_i)e_i
$$
et $[\nabla f(x), x] = 0$ s'\'ecrit : 
$$
	\sum_i\ B(\nabla f(x), e^*_i) [e_i,x] = 0 \quad \hbox{pour tout}\ x.
$$
On note $L_i : \mathfrak{g} \longrightarrow \mathfrak{g}$ le champ de vecteurs adjoint, $L_i(x) = [e_i,x]$
pour tout $x$. On a donc :
$$
	\sum_i (\partial_i f)L_i = 0
$$
o\`u $\partial_i f$ est la d\'eriv\'ee de $f$ le long du vecteur $e^*_i$. Pour de plus amples d\'etails et des
compl\'ements, le lecteur pourra se reporter \`a \cite{L-S}(lemma 2.2).

\vskip 10mm
\begin{enumerate}
\item Dans la suite, $\mathfrak{g} = \mathfrak{g}\ell(n,\Bbb{R})$, avec : 
$$
	B(x,y) = tr(xy)\quad (x,y) \in \mathfrak{g} \times \mathfrak{g}
$$
et $\mathfrak{p}= Lie(P)$. On utilise la base naturelle $(E_{ij})_{i,j}$ de $\mathfrak{g}$ et sa base duale
$(E^*_{ij})_{i,j}$, avec $E^*_{ij} = E_{ji}$, et on note $L_{ij}$ le champ adjoint : $L_{ij}(x) = [E_{ij},x]$. On a
donc : 
$$
	\sum_{i,j}\ tr(\nabla f(x) E_{ji})[E_{ij},x] = 0 \quad \hbox{pour tout}\ x
$$
lorsque $f$ est invariante.

	On applique ceci successivement aux fonctions invariantes : $p_k(x) = {1\over k}\, tr(x^k)\ (1\leq k \leq
n)$ de sorte que $\nabla p_k(x) = x^{k-1}$, et on obtient les $n$ \'egalit\'es : 
$$
	\sum_{i,j}\ tr(x^k E_{ji})L_{ij}(x) = 0\quad (0 \leq k \leq n-1).
$$
Soit $\varphi : \mathfrak{g} \longrightarrow \Bbb{R}$ une fonction de classe $C^1$. On a alors : 
$$
	\sum_{i,j}\, tr(x^k E_{ji})(L_{ij}\varphi)(x) = 0 \quad (0\leq k \leq n-1)
$$
($L_{ij}$ est consid\'er\'e comme un op\'erateur diff\'erentiel lin\'eaire homog\`ene de degr\'e 1). En
particulier, lorsque $\varphi$ est localement $P$-invariante, i.e. lorsque :
$$
	L_{ij} \varphi = 0\quad 1 \leq i \leq n-1,\ 1 \leq j \leq n,
$$
il reste : 
$$
	\sum_{1\leq j \leq n} tr(x^k\, E_{jn})L_{nj}\, \varphi(x) = 0 \quad (0\leq k \leq n-1).
$$
Il s'agit l\`a d'un syst\`eme lin\'eaire \`a $n$ inconnues $L_{nj}\, \varphi(x)$ (pour $x$ fix\'e) dont le
d\'eterminant est : 
$$
	D(x) = \hbox{d\'et}\big(tr(x^k E_{jn}) _{\begin{subarray}{1} 0 \leq k \leq n-1\\ 1\leq j \leq n
\end{subarray}}\big).
$$
La fonction $D$ est un polyn\^ome homog\`ene de degr\'e $n(n-1)/2$, qui est \textsl{non nulle}. On
remarque en effet que :
$$
	D(x) = [e_n, e_n \, x,\ldots , e_n \, x^{n-1}]
$$
est le d\'eterminant des $n$ vecteurs lignes $e_n, e_n \, x,\ldots , e_n \, x^{n-1}$. On constate alors que
lorsque $x=x_0$ est une ``matrice compagnon'' : 
$$
	x_0 =  \begin{array}{|lllll|}
0	&0	&\ldots	&\ldots	&\alpha_n\\
1	&0	&\ldots	&\ldots	&\alpha_{n-1}\\
0	&1	&0	&\ldots	&\alpha_{n-2}\\
	&	&\ddots		&\ddots &\vdots	\\
0	&0	&			&1		&\alpha_1\\
\end{array}
$$
on a : $D(x) = [e_n, e_{n-1},\ldots , e_1] = \pm 1$.

\smallskip
Il en r\'esulte que $L_{n1}\, \varphi(x) = L_{n2}\, \varphi(x) =\cdots = L_{nn}\, \varphi(x) = 0$ pour tout
$x$ tel que : $D(x)\not= 0$. On a donc :

\vskip 7mm
\noindent
\textbf{Lemme}~:~Toute fonction $\varphi$, de classe $C^1$ et localement $P$-invariante, est localement
$GL(n,\Bbb{R})$-invariante.

\medskip
\item  On note $\Omega$ l'ensemble des $x$ tels que $D(x) \not= 0$, c'est-\`a-dire l'ensemble des $x$ tels
que les vecteurs lignes $e_n, e_nx,\ldots ,e_nx^{n-1}$ soient lin\'eairement ind\'ependants. Donc :
$\Omega \subset \mathfrak{g}_r$, o\`u $\mathfrak{g}_r$ est l'ensemble des \'el\'ements r\'eguliers de
l'alg\`ebre de Lie $\mathfrak{g} = \mathfrak{g}\ell(n,\Bbb{R})$. Par ailleurs, si $y \in P$, on a : 
\begin{eqnarray}
D(yxy^{-1}) &=		&[e_n, e_nyxy^{-1},\ldots , e_nyx^{n-1}y^{-1}]\nonumber\\
	&=		&[e_n y^{-1}, e_nyxy^{-1},\ldots , e_nyx^{n-1}y^{-1}]\nonumber\\
&=	&(det\, y^{-1})[e_n, e_nyx,\ldots , e_nyx^{n-1}]\nonumber\\
&=	&(det \,y^{-1})[e_n, e_n x,\ldots , e_n x^{n-1}]
\nonumber
\end{eqnarray}
Donc $\Omega$ est un ouvert de Zariski, \textsl{$P$-invariant}, constitu\'e d'\'el\'ements r\'eguliers.

\medskip
\item  Soit $T$ une distribution localement $P$-invariante sur $\mathfrak{g}\ell(n,\Bbb{R})$, i.e. telle que :
$$
	L_{ij} T = 0 \quad \hbox{lorsque}\quad 1 \leq i \leq n-1,\quad 1\leq j \leq n.
$$
Comme ci-dessus, il vient :
$$
	\sum^n_{j=1}\, tr(x^k\ E_{jn})L_{nj} T = 0 \quad (0 \leq k \leq n-1)
$$
et par cons\'equent : 
$$
	D(L_{nj}T) = 0 \quad (1\leq j \leq n)
$$
o\`u, dans le premier membre de l'\'egalit\'e pr\'ec\'edente, figure le produit de la fonction $(C^\infty)D$
par la distribution $L_{nj}T$. Donc $L_{nj}T$ est nulle dans $\Omega$, pour tout entier $j$ v\'erifiant
$1\leq j \leq n$.

\vskip 7mm
\noindent
\textbf{Lemme}~:~Soit $T$ une distribution localement $P$-invariante sur
$\mathfrak{g}\ell(n,\mathbb{R})$. Alors $T$ est localement $\mathfrak{g}\ell(n,\mathbb{R})$-invariante
sur $\Omega$.

\medskip
\item  \textbf{Remarques} : 1) Le couple $(e_n,x)$ d\'efinit une forme lin\'eaire sur l'alg\`ebre de Lie du
groupe affine Aff($\mathbb{R}^n)$ de $\mathbb{R}^n$. La condition :
$$
	[e_n, e_nx,\ldots,e_nx^{n-1}] \not = 0
$$
exprime que cette forme lin\'eaire appartient \`a une orbite coadjointe ouverte.

\medskip
2) Comme remarqu\'e plus haut, l'ouvert $\Omega$ contient les ``matrices compagnons''. Par suite, le
satur\'e de $\Omega$ sous l'action adjointe de $GL(n,\mathbb{R})$ est l'ouvert $\mathfrak{g}_r$ des
\'el\'ements r\'eguliers de $\mathfrak{g}\ell(n,\mathbb{R})$.
\end{enumerate}

\vskip 12mm

\renewcommand{\refname}{Bibliographie} 

\end {document}